\documentclass[10pt, a4paper]{article}

\usepackage{amsmath,  amsfonts}

\usepackage[T1]{fontenc}

\newtheorem{theorem}{\quad Theorem}[section]

\newcommand{\be} {\begin{equation}}
\newcommand{\ee} {\end{equation}}

\title {Bounds for the volume of the solutions to a system on the annulus.}

\date{}

 \author{Samy Skander Bahoura\footnote {e-mails: samybahoura@yahoo.fr, samybahoura@gmail.com} \\ 
 {\small Equipe d'Analyse Complexe et G\'eom\'etrie.}\\  
  {\small Universit\'e Pierre et Marie Curie, 75005 Paris, France.}} 

\begin{document}

\maketitle
\begin{abstract}

We consider an elliptic system with regular H\"olderian weight and exponential nonlinearity or with weight and boundary singularity, and, Dirichlet condition. We prove the boundedness of the volume of the solutions to those systems on the annulus.

\end{abstract}



\section{Introduction and Main Results} 

We set $ \Delta = \partial_{11} + \partial_{22} $  on the annulus $ \Omega = C(1,1/2,0) $ of $ {\mathbb R}^2 $ of radii $ 1 $ and $ 1/2 $ centered at the origin.

\smallskip

We consider the following system:

$$ (P)   \left \{ \begin {split} 
       -\Delta u  & = (1+|x-x_0|^{2\beta}) V e^{v} \,\, &\text{in} \,\, & \Omega  \subset {\mathbb R}^2, \\
      - \Delta v  & = W e^{u} \,\, &\text{in} \,\, & \Omega  \subset {\mathbb R}^2, \\
                    u & = 0  \,\,             & \text{on} \,\,    &\partial \Omega,\\ 
                    v & = 0  \,\,             & \text{on} \,\,    &\partial \Omega.                        
   \end {split}\right.
$$

Here: $ C(1) $ the unit circle and $ C(1/2) $ the circle of radius $ 1/2 $ centered at the origin.

$$ \beta \in (0, 1/2), \,\, x_0 \in C(1/2), $$
and,

$$ u \in W_0^{1,1}(\Omega), \,\, e^u \in L^1({\Omega}) \,\, {\rm and} \,\,  0 < a \leq V \leq b, $$

and,

$$ v \in W_0^{1,1}(\Omega), \,\, e^v \in L^1({\Omega}) \,\, {\rm and} \,\,  0 < c \leq W \leq d. $$

This is a system with regular H\"olderian weight not Lipschitz in $ x_0 $ but have a weak derivative.

This problem $ (P) $ is defined in the sense of the distributions, see [10]. The system was studied by many authors, see [14, 16, 27], also for Riemannian surfaces, see [1-27],  where one can find some existence and compactness results. In [9] we have an interior estimate for elliptic equations with exponential nonlinearity.

In this paper we try to prove that we have on all $ \Omega $ the boundedness of the volume of the solutions of $ (P) $ if we add the assumption that $ V $ and $ W $ are uniformly Lipschitz with particular Lispchitz numbers.

Here we have:

\begin{theorem} Assume that $ u $ is a solution of $ (P) $ relative to $ V $ and $ W $ with the following conditions:

$$  x_0 \in C(1/2) \subset \partial \Omega, \,\, \beta \in (0,1/2), $$

and,

$$ 0 < a \leq V \leq b,\,\, ||\nabla V||_{L^{\infty}} \leq A=\dfrac{a}{2(1+2^{2\beta})}, $$

and,

$$ 0 < c \leq W \leq d,\,\, ||\nabla W||_{L^{\infty}} \leq B=\dfrac{c}{2}, $$

we have,

$$  \int_{\Omega} e^u \leq c(a, b, c, d, \beta, x_0, \Omega), \,\, {\rm and}, \,\, \int_{\Omega} e^v \leq c'(a, b, c, d, \beta, x_0, \Omega)  $$

\end{theorem} 

We have the same result if we consider a system with boundary singularity. On the annulus $ \Omega = C(1,1/2,0) $ of $ {\mathbb R}^2 $ of radii $ 1 $ and $ 1/2 $ centered at the origin.

We consider the following system:

$$ (P_{\beta})   \left \{ \begin {split} 
       -\Delta u  & = |x-x_0|^{2\beta} V e^{v} \,\, &\text{in} \,\, & \Omega  \subset {\mathbb R}^2, \\
      - \Delta v  & = |x-x_0|^{2\beta} W e^{u} \,\, &\text{in} \,\, & \Omega  \subset {\mathbb R}^2, \\
                    u & = 0  \,\,             & \text{on} \,\,    &\partial \Omega,\\ 
                    v & = 0  \,\,             & \text{on} \,\,    &\partial \Omega.                        
   \end {split}\right.
$$

Here: $ C(1) $ the unit circle and $ C(1/2) $ the circle of radius $ 1/2 $ centered at the origin.

$$ \beta \in (-1/2, +\infty), \,\, x_0 \in C(1/2), $$

and,

$$ u \in W_0^{1,1}(\Omega), \,\, |x-x_0|^{2\beta} e^u \in L^1({\Omega}) \,\, {\rm and} \,\,  0 < a \leq V \leq b, $$

and,

$$ v \in W_0^{1,1}(\Omega), \,\, |x-x_0|^{2\beta} e^v \in L^1({\Omega}) \,\, {\rm and} \,\,  0 < c \leq W \leq d. $$

Here we have:

\begin{theorem} Assume that $ u $ is a solution of $ (P_{\beta}) $ relative to $ V $ and $ W $ with the following conditions:

$$  x_0 \in C(1/2) \subset \partial \Omega, \,\, \beta \in (-1/2, +\infty), $$

and,

$$ 0 < a \leq V \leq b,\,\, ||\nabla V||_{L^{\infty}} \leq A=\dfrac{(\beta+1)a}{2}, $$

and,

$$ 0 < c \leq W \leq d,\,\, ||\nabla W||_{L^{\infty}} \leq B=\dfrac{(\beta+1)c}{2}, $$

we have,

$$  \int_{\Omega} |x-x_0|^{2\beta} e^u \leq c(a, b, c, d, \beta, x_0, \Omega), \,\, {\rm and}, \,\, \int_{\Omega} |x-x_0|^{2\beta} e^v \leq c'(a, b, c, d, \beta, x_0, \Omega)  $$

\end{theorem} 

We have the same result if we consider a system with boundary singularity. On the annulus $ \Omega = C(1,1/2,0) $ of $ {\mathbb R}^2 $ of radii $ 1 $ and $ 1/2 $ centered at the origin.

We consider the following system:

$$ (P_{\beta})   \left \{ \begin {split} 
       -\Delta u  & = |x-x_0|^{2\beta} V e^{v} \,\, &\text{in} \,\, & \Omega  \subset {\mathbb R}^2, \\
      - \Delta v  & = W e^{u} \,\, &\text{in} \,\, & \Omega  \subset {\mathbb R}^2, \\
                    u & = 0  \,\,             & \text{on} \,\,    &\partial \Omega,\\ 
                    v & = 0  \,\,             & \text{on} \,\,    &\partial \Omega.                        
   \end {split}\right.
$$

Here: $ C(1) $ the unit circle and $ C(1/2) $ the circle of radius $ 1/2 $ centered at the origin.

$$ \beta \in (-1/2, +\infty), \,\, x_0 \in C(1/2), $$

and,

$$ u \in W_0^{1,1}(\Omega), \,\,  e^u \in L^1({\Omega}) \,\, {\rm and} \,\,  0 < a \leq V \leq b, $$

and,

$$ v \in W_0^{1,1}(\Omega), \,\, |x-x_0|^{2\beta} e^v \in L^1({\Omega}) \,\, {\rm and} \,\,  0 < c \leq W \leq d. $$

Here we have:

\begin{theorem} Assume that $ u $ is a solution of $ (P_{\beta}) $ relative to $ V $ and $ W $ with the following conditions:

$$  x_0 \in C(1/2) \subset \partial \Omega, \,\, \beta \in (-1/2, +\infty), $$

and,

$$ 0 < a \leq V \leq b,\,\, ||\nabla V||_{L^{\infty}} \leq A=\dfrac{(\beta+1)a}{2}, $$

and,

$$ 0 < c \leq W \leq d,\,\, ||\nabla W||_{L^{\infty}} \leq B=\dfrac{c}{2}, $$

we have,

$$  \int_{\Omega} e^u \leq c(a, b, c, d, \beta, x_0, \Omega), \,\, {\rm and}, \,\, \int_{\Omega} |x-x_0|^{2\beta} e^v \leq c'(a, b, c, d, \beta, x_0, \Omega)  $$

\end{theorem} 

\section{Proof of the Theorems:}

\smallskip

{\it Proof of the theorem 1.1:}

\smallskip

By corollary 1 of the paper of Brezis-Merle, we have: $ e^{ku}, e^{kv} \in L^1(\Omega) $ for all $ k >2 $ and the elliptic estimates and the Sobolev embedding imply that: $ u,v \in W^{2,k}(\Omega) \cap C^{1,\epsilon}(\bar \Omega), \epsilon >0 $. By the maximum principle $ u,v \geq 0 $.

\smallskip

Step 1: We use the first eigenvalue and the first eigenfunction with Dirichlet boundary condition to bound the volumes locally uniformly. Thus the solutions are locally uniformly bounded by Brezis-Merle result. The solutions $ u,v >0 $ are locally uniformly bounded in $ C^{1,\epsilon}(\Omega) $ for $ \epsilon $ small.

\smallskip

By Cauchy-Schwarz inequality, applied to $ u \sqrt {\phi_1} $ and $ {\sqrt \phi_1} $ for the following equality:

$$ \int_{\Omega} (1+|x-x_0|^{2\beta}) Ve^v \phi_1dx =\lambda_1 \int_{\Omega} u\phi_1 dx \leq c_1 {(\int_{\Omega} u^2 \phi_1)}^{1/2} \leq c_2 {(\int_{\Omega} W e^u \phi_1)}^{1/2}, $$

and for $ v $,

$$ \int_{\Omega} We^u \phi_1dx =\lambda_1 \int_{\Omega} v\phi_1 dx \leq c_3 {(\int_{\Omega} v^2 \phi_1)}^{1/2} \leq c_4 {(\int_{\Omega} (1+|x-x_0|^{2\beta })V e^v \phi_1)}^{1/2} , $$

Thus,

$$ (\int_{\Omega} (1+|x-x_0|^{2\beta})Ve^v \phi_1dx)^{3/4} \leq c_5, $$

and,

$$ (\int_{\Omega} We^u \phi_1dx)^{3/4} \leq c_6, $$

We can use Brezis-Merle arguments to prove that for all subdomain $ K $ of $ \Omega $: the two integrals, $ \int_{K} (1+|x-x_0|^{2\beta})Ve^v dx, $ and, $ \int_{K} We^u dx $ converge to nonegative measues $ \mu_1, \mu_2 $ without nonregular points.

\smallskip

By contradiction, suppose that $ \max_K u_i \to +\infty $ and $ \max_K v_i \to +\infty $.

\smallskip

Since $ (1+|x-x_0|^{2\beta})V_ie^{v_i} $ and $ W_ie^{u_i} $ are bounded in $ L^1(K) $, we can extract from those two sequences two subsequences which converge to two nonegative measures $ \mu_1 $ and $ \mu_2 $. (This procedure is similar to the procedure of Brezis-Merle, we apply corollary 4 of Brezis-Merle paper, see [9]).

If $ \mu_1(y_0) < 4 \pi $, by a Brezis-Merle estimate for the first equation, we have $ e^{u_i} \in L^{1+\epsilon} $ uniformly around $ y_0 $, by the elliptic estimates, for the second equation, we have $ v_i \in W^{2, 1+\epsilon} \subset L^{\infty} $ uniformly around $ y_0 $, and , returning to the first equation, we have $ u_i \in L^{\infty} $ uniformly around $ y_0 $.

If $ \mu_2(y_0) < 4 \pi $, then $ u_i $ and $ v_i $ are also locally uniformly bounded around $ y_0 $.

Thus, we take a look to the case when, $ \mu_1(y_0) \geq 4 \pi $ and $ \mu_2(y_0) \geq 4 \pi $. By our hypothesis, those points $ y_0 $ are finite.

We will see that inside $ K $ no such points exist. By contradiction, assume that, we have $ \mu_1(y_0) \geq 4 \pi $. Let us consider a ball $ B_R(y_0) $ which contain only $ y_0 $ as nonregular point. Thus, on $ \partial B_R(y_0) $, the two sequence $ u_i $ and $ v_i $ are uniformly bounded. Let us consider:

$$ \left \{ \begin {split} 
       -\Delta z_i & = (1+|x-x_0|^{2\beta})V_i e^{v_i} \,\, &\text{in} \,\, & B_R(y_0)  \subset {\mathbb R}^2, \\
                       z_i & = 0  \,\,             & \text{in} \,\,    &\partial B_R(y_0).                                            
   \end {split}\right.
$$

By the maximum principle:

$$ z_i \leq u_i, $$ 

and $ z_i \to z $ almost everywhere on this ball, and thus,

$$ \int e^{z_i} \leq \int e^{u_i} \leq C, $$

and,

$$ \int e^z \leq C.$$

but, $ z  $ is  a solution in $ W_0^{1,q}(B_R(y_0)) $, $ 1\leq q <2 $, of the following equation:

$$ \left \{ \begin {split} 
       -\Delta z & = \mu_1\,\, &\text{in} \,\, & B_R(y_0) \subset {\mathbb R}^2, \\
                     z & = 0  \,\,             & \text{in} \,\,    &\partial B_R(y_0).                                            
   \end {split}\right.
$$

with, $ \mu_1 \geq 4 \pi $ and thus, $ \mu_1 \geq 4\pi \delta_{y_0} $ and then, by the maximum principle in $ W_0^{1,1}(B_R(y_0)) $:

$$ z \geq -2 \log |x-y_0|+ C $$

thus,

$$ \int e^z = + \infty, $$

which is a contradiction. Thus, there is no nonregular points inside  $ K $. Thus $(u_i)$ and $(v_i)$ are uniformly bounded in $ K $ and also in $ C^{1,\epsilon}(K) $ by the elliptic estimates, for all $ K \subset \subset \Omega $.

\smallskip

Step 2: Let's consider $ C_1=C(1,3/4,0) $ and $ C_2=C(3/4,1/2,0) $ the two annulus wich are the neighborhood of the two components of the boundary.

We multiply the equation by $ (x-x_0)\cdot \nabla u $ on $ C_1 $ and $ C_2 $ and use the Pohozaev-Rellich identity and Stokes theorem, see [26]. We use the fact that $ u $ and $ v $ are uniformly bounded around the circle $ C(3/4) $. We obtain:

1) We have on $ C_1 $:

$$ \int_{C_1} (\Delta u)[(x-x_0)\cdot \nabla v] dx=\int_{C_1} -[(1+|x-x_0|^{2\beta}) V (x-x_0)\cdot \nabla (e^v)]dx, $$

and,

$$ \int_{C_1} (\Delta v)[(x-x_0)\cdot \nabla u] dx=\int_{C_1} -[W (x-x_0)\cdot \nabla (e^u)]dx, $$

Thus, by integration by parts,

$$ \int_{C_1} (\Delta u)[(x-x_0)\cdot \nabla v]+ (\Delta v)[(x-x_0)\cdot \nabla u] dx = $$

$$
 \int_{\partial C_1} [(x-x_0)\cdot \nabla u ] (\nabla v \cdot \nu)+[(x-x_0)\cdot \nabla v ] (\nabla u \cdot \nu) -[(x-x_0)\cdot \nu] (\nabla u \cdot \nabla v)= $$
 
$$ = \int_{C_1} (2+2(\beta + 1)|x-x_0|^{2\beta}) Ve^v dx + \int_{C_1} (1+|x-x_0|^{2\beta}) (x-x_0)\cdot \nabla V e^v dx+  $$
 
$$ + \int_{C_1} 2We^u dx + \int_{C_1} (x-x_0)\cdot \nabla W e^u dx+  $$

$$ - \int_{\partial C_1} (1+|x-x_0|^{2\beta}) [(x-x_0)\cdot \nu ] V e^v d\sigma +$$

$$ - \int_{\partial C_1} [(x-x_0)\cdot \nu ] W e^u d\sigma $$

We can write, ($ u= 0 $ on $ C(1) $):

$$ \int_{C(1)} [(x-x_0)\cdot \nu ] (\partial_{\nu} u)(\partial_{\nu} v) d\sigma + O(1)= $$

$$ \leq k_1 \int_{C_1} (1+|x-x_0|^{2\beta}) Ve^v dx +k_2 \int_{C_1} We^u dx+ O(1) = $$

$$ = k_1\int_{C(1)} \partial_{\nu} u d\sigma + k_2\int_{C(1)} \partial_{\nu} v d\sigma +O(1), $$

with $ k_1,k_2>0 $ not depends on $ u $.

Be cause $ \nu = x, ||x||=1, ||x_0||=1/2 $ and then by Cauchy-Schwarz, $ (x-x_0)\cdot x=||x||^2-x_0 \cdot x \geq 1/2 $, we obtain:

\be 0 < \int_{C(1)} (\partial_{\nu} u ) (\partial_{\nu} v) d\sigma \leq  k_1\int_{C(1)} \partial_{\nu} u d\sigma + k_2\int_{C(1)} \partial_{\nu} v d\sigma +O(1), \ee

Let $ \epsilon =\inf (1, \frac{a}{d}) $, then: $ -\Delta (\epsilon v)=\epsilon W e^u \leq a e^u \leq (1+|x-x_0|^{2\beta})V e^u $, and $ \epsilon \leq  1, v\geq 0 \Rightarrow (1+|x-x_0|^{2\beta} ) V e^v \geq (1+|x-x_0|^{2\beta})V e^{\epsilon v} $ and then: $ -\Delta u \geq (1+|x-x_0|^{2\beta})V e^{\epsilon v} $. (We can remove $ (1+|x-x_0|^{2\beta})$ ) We obtain:

$$ -\Delta (u-\epsilon v)= -\Delta u+ \Delta (\epsilon v) \geq (1+|x-x_0|^{2\beta})V(e^{\epsilon v}-e^u), $$

Thus,

$$ -\Delta (u-\epsilon v) + (1+|x-x_0|^{2\beta})V(e^u-e^{\epsilon v}) \geq 0, $$

(For the theorem 1.2, we have the weight $ |x-x_0|^{2\beta} $, in the two equations, we can compare $ \epsilon v $ and $ u $).

We return to the proof of theorem 1.1. Let's consider the fonction:

$$ c(x)=(1+|x-x_0|^{2\beta})V\dfrac{(e^u-e^{\epsilon v})}{u-\epsilon v}, \, {\rm if} \, u\not = \epsilon v, \, and, \, c(x)=(1+|x-x_0|^{2\beta})Ve^u,\, {\rm if} \, u=\epsilon v. $$

The function $ c \geq 0 $ is $ C^{\beta}(\bar \Omega) $, and $ -c \leq 0 $.(For the theorem 1.2, if $ -1/2< \beta <0$, $ c $ is $ C^{-\beta}(\Omega)\cap L^2(\Omega) $, but this is sufficient to apply the weak maximum principle, see the book of Gilbarg-Trudinger).

We can write:

$$ \Delta (\epsilon v-u)- c(x) (\epsilon v-u) \geq 0, \, {\rm in } \, \Omega, \, and, \, u-\epsilon v=0 \, {\rm on} \,\, \partial \Omega, $$

The operator $ L= \Delta + (-c) =\Delta + \tilde c $ satisfies the maximum principle because $ \tilde c =-c \leq 0 $, we obtain:

by the weak maximum principle for $ \epsilon v-u \in C^2(\Omega) \cap C^1(\bar \Omega) $, see the book of Gilbarg-Trudinger, we obtain (and for the outer normal):

$$ \epsilon v - u \leq 0, \, {\rm and \, then }, \, \partial_{\nu} (\epsilon v-u) \geq 0, $$

Thus, for the inner normal, we have:

\be \partial_{\nu} u \geq \epsilon \cdot \partial_{\nu} v >0 \,\, {\rm on } \,\, \partial \Omega. \ee

By the same argument, if we set $ \bar \epsilon= \inf (1, \frac{c}{(1+2^{2\beta})b}) $, we obtain, for the inner normal:

\be \partial_{\nu} v \geq \bar \epsilon \cdot \partial_{\nu} u >0 \,\, {\rm on } \,\, \partial \Omega. \ee

{\bf Remark:} For theorems 1.2 and 1.3: we can remove the fact that we have the weight $ |x-x_0|^{2\beta} $ in the two equations of the system, we can assume that we have one equation with weight and the other equation without weight for example. Indeed, remark that $ x_0 \in C(1/2) $ and here we consider $ C(1) $, we can apply the weak maximum principle on $ C_1 $ after choosing $ 1> \epsilon >0 $ and $ 1 > \bar \epsilon >0 $ independent of $ u $ and $ v $, such that: $ \epsilon v - u \leq 0 $ on $ C(3/4) $ and $ \bar \epsilon u-v \leq 0 $ on $ C(3/4) $, because we have the uniform interior estimate around $ C(3/4) $ and the maximum principle (by contradiction if we assume there are $ (t_i) \in C(3/4), t_i \to t_0 $ and $ u_i, v_i $ such that $ u_i(t_i) \to 0 $, $ (u_i), (v_i) $ converge to $ u_0,v_0 $ on $ B(t_0, r_0), r_0 >0 $ with $ u_0(t_0)=0 $, but $ -\Delta u_0 >0 $ and the maximum principle applied to $ -u_0 $ imply that $ u_0 >0 $ around $ t_0 $. This fact imply that $ u, v $ have positive lower bounds on $ C(3/4)$ and we can choose $ \epsilon, \bar \epsilon $ independent of $ u, v $ on $C(3/4)$. Finally we can apply the weak maximum principle on $ C_1 $).

Thus, we use $ (1), (2), (3) $, we obtain the same result as for one equation;

$$ \int_{C(1)} (\partial_{\nu} v)^2 d\sigma \leq C_3 \int_{C(1)} \partial_{\nu} v d\sigma + C_4, $$

with, $ C_3, C_4 >0 $ and not dpend on $ u $ and $ v $, and,

$$ \int_{C(1)} (\partial_{\nu} u)^2 d\sigma \leq C_5 \int_{C(1)} \partial_{\nu} u d\sigma + C_6, $$

with, $ C_5, C_6 >0 $ and not dpend on $ u $ and $ v $, and,

Applying the Cauchy-Schwarz inequality, we have:

$$  \int_{C(1)} (\partial_{\nu} u)^2 d\sigma =O(1), \,\, \int_{C(1)} (\partial_{\nu} v)^2 d\sigma = O(1) $$

and thus,

$$  \int_{C(1)} \partial_{\nu} u d\sigma =O(1), \,\, \int_{C(1)} \partial_{\nu} v d\sigma = O(1) $$

2) We have on $ C_2 $, we use again, the uniform boundedness of $ u $ and $ v $ in $ C^1 $ norm around $ C(3/4) $:

$$ \int_{C_2} (\Delta u)[(x-x_0)\cdot \nabla v] dx=\int_{C_2} -[(1+|x-x_0|^{2\beta}) V (x-x_0)\cdot \nabla (e^v)]dx, $$

and,

$$ \int_{C_2} (\Delta v)[(x-x_0)\cdot \nabla u] dx=\int_{C_2} -[W (x-x_0)\cdot \nabla (e^u)]dx, $$

Thus, by integration by parts,

$$ \int_{C_2} (\Delta u)[(x-x_0)\cdot \nabla v]+ (\Delta v)[(x-x_0)\cdot \nabla u] dx = $$

$$
 \int_{\partial C_2} [(x-x_0)\cdot \nabla u ] (\nabla v \cdot \nu)+[(x-x_0)\cdot \nabla v ] (\nabla u \cdot \nu) -[(x-x_0)\cdot \nu] (\nabla u \cdot \nabla v)= $$
 
$$ = \int_{C_2} (2+2(\beta + 1)|x-x_0|^{2\beta}) Ve^v dx + \int_{C_2} (1+|x-x_0|^{2\beta}) (x-x_0)\cdot \nabla V e^v dx+  $$
 
$$ + \int_{C_2} 2We^u dx + \int_{C_2} (x-x_0)\cdot \nabla W e^u dx+  $$

$$ - \int_{\partial C_2} (1+|x-x_0|^{2\beta}) [(x-x_0)\cdot \nu ] V e^v d\sigma +$$

$$ - \int_{\partial C_2} [(x-x_0)\cdot \nu ] W e^u d\sigma $$

But here, $ \nu=-2x$ and $ (x-x_0)\cdot \nu= -2(x-x_0)\cdot x \leq 0 $ and thus:

$$  \int_{C(1/2)} (x-x_0)\cdot \nu (\partial_{\nu} u ) (\partial_{\nu} v) d\sigma + O(1) = $$

$$ = \int_{C_2} (2+2(\beta + 1)|x-x_0|^{2\beta}) Ve^v dx + \int_{C_2} (1+|x-x_0|^{2\beta}) (x-x_0)\cdot \nabla V e^v dx+  $$
 
$$ + \int_{C_2} 2We^u dx + \int_{C_2} (x-x_0)\cdot \nabla W e^u dx+  $$

$$ - \int_{\partial C_2} (1+|x-x_0|^{2\beta}) [(x-x_0)\cdot \nu ] V e^v d\sigma +$$

$$ - \int_{\partial C_2} [(x-x_0)\cdot \nu ] W e^u d\sigma $$

thus:

$$ \int_{C_2} (2+2(\beta +1)|x-x_0|^{2\beta}) Ve^v dx + \int_{C_2} (1+|x-x_0|^{2\beta}) (x-x_0)\cdot \nabla V e^v dx + $$

$$ + \int_{C_2} 2We^u dx + \int_{C_1} (x-x_0)\cdot \nabla W e^u dx+  $$

$$ + \int_{C(1/2)} \frac{1}{2}[-(x-x_0)\cdot \nu ] (\partial_{\nu} u)(\partial_{\nu} v) d\sigma = O(1) $$

If we choose:

\be \dfrac{|(x-x_0)\cdot \nabla V|}{V} \leq \dfrac{1}{2} \inf_{x \in \bar \Omega } \dfrac{(2+2(\beta+1)|x-x_0|^{2\beta})}{1+|x-x_0|^{2\beta}}, \ee

and,

\be \dfrac{|(x-x_0)\cdot \nabla W|}{W} \leq \dfrac{1}{2} \cdot 2, \ee

(For the theorem 1.2, we choose: $ \dfrac{|(x-x_0)\cdot \nabla V|}{V} \leq \dfrac{1}{2} \inf_{x \in \Omega } \dfrac{(2(\beta+1)|x-x_0|^{2\beta})}{|x-x_0|^{2\beta}}=\beta+1,$ and also for $ W $, we choose $ \dfrac{|(x-x_0)\cdot \nabla W|}{W} \leq \dfrac{1}{2} \inf_{x \in \Omega } \dfrac{(2(\beta+1)|x-x_0|^{2\beta})}{|x-x_0|^{2\beta}} =(\beta+1) $).

We obtain:

$$ \int_{C_2} (2+2(\beta +1)|x-x_0|^{2\beta}) Ve^v dx + \int_{C_2} (1+|x-x_0|^{2\beta}) (x-x_0)\cdot \nabla V e^v dx \geq $$

$$ \geq \frac{1}{2} \int_{C_2} (2+2(\beta +1)|x-x_0|^{2\beta}) Ve^v dx \geq 0 $$

and,

$$ \int_{C_2} 2 W e^u dx + \int_{C_2} (x-x_0)\cdot \nabla W e^u dx \geq $$

$$ \geq \frac{1}{2} \int_{C_2} 2 W e^u dx \geq 0 $$

thus,

$$ \int_{C_2} [(2+2(\beta +1)|x-x_0|^{2\beta}) Ve^v+ 2W e^u] dx =O(1), $$

we obtain:
$$ \int_{C_2} (1+|x-x_0|^{2\beta}) Ve^v dx =O(1), $$

and,

$$ \int_{C_2} W e^u dx =O(1), $$

and thus, if $  A\leq \dfrac{a}{2(1+2^{2\beta})} $ and $ B \leq \dfrac{c}{2} $, we obtain :

$$ \int_{C(1/2)} \partial_{\nu} u d\sigma = O(1), \,\,  \int_{C(1/2)} \partial_{\nu} v d\sigma = O(1), $$

Thus, if we use 1) and 2), we obtain: if $ A\leq \dfrac{a}{2(1+2^{2\beta})} $ and $ B \leq \dfrac{c}{2} $:

$$\int_{C(1)} (\partial_{\nu} u) d\sigma =O(1),\,\, {\rm and} \,\, \int_{C(1/2)} \partial_{\nu} u d\sigma = O(1), $$

and,

$$\int_{C(1)} (\partial_{\nu} v) d\sigma =O(1),\,\, {\rm and} \,\, \int_{C(1/2)} \partial_{\nu} v d\sigma = O(1), $$

and, thus:

$$ \int_{\Omega } [(1+|x-x_0|^{2\beta}) Ve^v] dx=\int_{\partial \Omega} (\partial_{\nu} u) d\sigma = O(1). $$

and,

$$ \int_{\Omega } We^u dx=\int_{\partial \Omega} (\partial_{\nu} v) d\sigma =O(1). $$

For the theorem 1.2. we obtain the same result if $ A\leq \dfrac{(\beta+1) a}{2} $ and $ B \leq \dfrac{(\beta+1)c}{2} $.

We have the same result for theorem 1.3.

\end{document}